\documentclass[11pt,reqno]{article}
\usepackage{bm}
\usepackage{amsthm}
\usepackage{amsmath}
\usepackage{amssymb}
\usepackage[round]{natbib}
\usepackage{color}
\usepackage{authblk}
\usepackage{microtype}
\usepackage{bm}
\usepackage{graphicx}
\usepackage{float}

\usepackage{geometry}
\newcommand{\rnc}{\renewcommand}
\newcommand{\nc}{\newcommand}
\newcommand{\mrm}{\mathrm}

\nc{\mb}{\mathbb}
\nc{\mc}{\mathcal}
\nc{\E}{\mb{E}}
\nc{\N}{\mb{N}}
\nc{\R}{\mb{R}}
\nc{\Q}{\mb{Q}}
\rnc{\P}{\mrm P}
\rnc{\d}{\mrm d}
\nc{\C}{\mc{C}}
\nc{\D}{\mc{D}}
\nc{\B}{\mc{B}}
\nc{\vbeta}{\bm \beta}
\nc{\vtheta}{\bm \theta}
\nc{\vX}{\bm X}
\nc{\vy}{\bm y}
\nc{\vU}{\bm U}
\nc{\vI}{\bm I}
\nc{\vE}{\bm E}
\nc{\ve}{\bm e}
\nc{\vV}{\bm V}
\nc{\vv}{\bm v}
\nc{\vS}{\bm S}
\nc{\vSigma}{\bm \Sigma}
\nc{\oPo}{\stackrel{\mrm p}{\rightarrow}}
\nc{\oWo}{\stackrel{w}{\rightarrow}}
\nc{\oDo}{\stackrel{d}{\longrightarrow}}
\nc{\eff}{\|F\|}

\def\E{{ E }}
\def\R{{ \mathbb{R} }}
\def\N{{ \mathbb{N} }}
\def\var{ \text{Var} }
\def\P{ P }
\def\ho{{ H_=} }
\def\kgeq{ K_{\geq}}

\def\E{ E }
\def\var{ \text{var} }

\newtheorem{theorem}{Theorem}

\newcommand\blfootnote[1]{%
  \begingroup
  \renewcommand\thefootnote{}\footnote{#1}%
  \addtocounter{footnote}{-1}%
  \endgroup
} 

\begin{document}

\title{\Large \bf 
Wild bootstrap logrank tests with broader power functions\\for testing superiority}
\author[1$*$]{Marc Ditzhaus}
\author[1]{Markus Pauly}

\affil[1]{Institute of Statistics, Ulm University, Germany.}

\maketitle

\begin{abstract}
\blfootnote{${}^*$ e-mail: marc.ditzhaus@hhu.de}
We introduce novel wild bootstrap procedures for testing superiority in unpaired two-sample survival data. By combining different classical weighted logrank test we obtain tests with broader power behavior. Right censoring within the data is allowed and may differ between the groups. The tests are shown to be asymptotically exact under the null, consistent for fixed alternatives and admissible for a larger set of local alternatives. 
Beside these asymptotic properties we also illustrate the procedures' strength in simulations for finite sample sizes. The tests are implemented in the novel \textsc{R}-package mdir.logrank and its application is demonstrated in an exemplary data analysis. 
\end{abstract}




\section{Introduction}\label{sec:intro}

A clinical two sample study may be conducted to infer superiority of a treatment over a control. For time to event data, this task is usually coped with the one-sided version of the logrank test \citep{Mantel1966, PetoPeto1972}. 
As it is 'only' (asymptotically) optimal to detect alternatives in which the hazards are proportional, several substitutes have been proposed. Beneath weighted logrank tests \citep{TaroneWare1977, Gill1980, HarringtonFleming1982, FlemingHarrington, lachin2009biostatistical, BagonaviciusETAL2010}, different modifications have been established in order to derive certain power 'enhancements' (at least in the two-sided case), see, e.g., \cite{Ehm1995, KosorokLin1999, BathkeETAL2009, YangPrentice2010, BrendelETAL2014, GaresETAL2015} and the references cited therein. From a methodological point of view, the projection-type procedure of \cite{BrendelETAL2014} slightly stands out as it in principle allows the simple construction of a permutation test that is (asymptotically) optimal against a pre-specified number of alternatives of interest. However, the method has not yet found its way into statistical practice for several reasons: 
\begin{itemize}
 \item[1.] It primarily focus on a mathematical sound development of the method; rigorously using functional analytic terminology. 
\end{itemize}
For example, the test statistic is derived as projection of a vector of weighted logrank statistics onto a functional cone corresponding to the choices of alternatives of interest. Although this is greatly appreciated by the current authors its abstract description and very complex (and non-pivotal) limit distribution might be too distant for most applied statisticians; especially 
\begin{enumerate}
 \item[2.] since the method has not yet been implemented in statistical software not to mention plain user manuals.
 \item[3.] Moreover, due to the complex limit distribution of their test statistic, the proposed permutation method utilizes a rather unusual studentization technique by writing the critical value as part of the test statistic to ensure asymptotic correctness.
\end{enumerate}
While it is nowadays (more or less) accepted that permutation tests for complex heterogeneous models (here given by possibly different censoring distributions) need a certain studentization \citep{Neuhaus1993, janssen1997studentized, JanssenMayer2001, chung2013exact, pauly2015asymptotic} the permutation technique of  \cite{BrendelETAL2014} typically needs two iterative Monte-Carlo steps: Beneath the usual Monte-Carlo approximation of the permutation distribution an additional one is a priori needed to calculate the critical value within the test statistic. This regrettably results in too time consuming calculations for critical values. 

It is the aim of the present paper to address these points accordingly. In particular, employing multiplier wild bootstrap resampling \citep{Lin1997, BeyersmannETAL2013, dobler2014bootstrapping, dobler2017non, bluhmki2018wild} instead of permutation enables us to directly cope with the complex limit distribution without making the above studentization detour. We rigorously analyze the asymptotic properties of the resulting procedure, specifically preserving the asymptotic optimality of the
\cite{BrendelETAL2014} method (Section~\ref{sec:wild} below). To give recommendations for choosing proper multipliers and analyze the procedures' finite sample properties, extensive simulations (Section~\ref{sec:simus}) are conducted under the null hypothesis and various alternatives. 
For an easy application, we implemented the resulting recommendations within the \textsc{R}-package {\bf mdir.logrank} available on GitHub. 
Therein, our novel procedure can either be called via a user-friendly GUI (based on default choices of typical alternatives) or a more enhanced function allowing for user-specific definitions of more specific alternatives. Its application is exemplified during an analysis of lung cancer data (Section~\ref{sec:realdata}). 
Before we follow this course of action we first fix the basic set-up in Section~\ref{sec:set-up}. 

We note that all proofs are deferred to the appendix, where also additional simulation results {and a demonstration of the graphical user interface} are presented.

\section{Two-sample survival set-up}\label{sec:set-up}
We consider a two-sample survival set-up given by mutually independent non-negative random variables 
\begin{align}\label{eqn:model}
	T_{j,i}\sim F_j \quad\text{ and }\quad C_{j,i}\sim G_{j} \quad (j=1,2;\,i=1,\ldots,n_j).
\end{align}
Here, $F_j$ and $G_j$ are continuous distribution functions, $T_{j,i}$ denotes the survival time of subject $i$ in group $j$ and 
$C_{j,i}$ its corresponding right-censoring time. As usual within this setting, only the right-censored survival times $X_{j,i}=\min\{T_{j,i},C_{j,i}\}$ and their corresponding censoring status $\delta_{j,i}=\mathbf{1}\{X_{j,i}=T_{j,i}\}$ $(j=1,2;\,i=1,\ldots,n_j)$  are actually observable.  Based on the resulting $n=n_1+n_2$ pooled observations we would now like to infer whether the second group (e.g., the intervention group) is superior to the first group (e.g., the control group) in terms of survival times. 
This leads to testing $\ho: \{A_1 = A_2\} : \{F_1 = F_2\}$ against the one-sided alternative of {\it stochastic ordering} or {\it superiority}
\begin{align}\label{eqn:testprob}
	\kgeq:\{A_1\geq A_2, A_1\neq A_2\} = \{F_1\geq F_2, F_1\neq F_2\},
\end{align}
where $A_j(t)=-\log(1-F_j(t))=\int_0^t(1-F_j)^{-1}\,\mathrm{ d }F_j$ denotes the cumulative hazard function, $j=1,2$, and 
the unknown censoring distributions can be interpreted as nonparametric nuisance parameters. In case of proportional hazards $dA_1(t) / dA_2(t) \equiv \lambda_1(t)/\lambda_2(t) \equiv \vartheta +1$ the well-known alternative $\{\vartheta>0\}$ particularly implies $\kgeq$.


To formulate relevant test statistics we adopt the usual counting process notation. Thus, let $N_{j,i}(t)=\mathbf{1}\{X_{j,i}\leq t,\,\delta_{j,i}=1\}$ and $Y_{j,i}(t)=\mathbf{1}\{X_{j,i}\geq t\}, \;$ $j=1,2;\,i=1,\ldots,n_j$. Summation over all subjects in group $j$ yields $N_j(t)=\sum_{i=1}^{n_j}N_{j,i}(t)$, the number of observed events in group $j$ until time $t$, and $Y_j(t)=\sum_{i=1}^{n_j}Y_{j,i}(t)$, the number of individuals in group $j$ being at risk just before $t$. 
The pooled processes $N=N_1+N_2$ and $Y=Y_1+Y_2$ can be interpreted in the same way. Now, the group-specific Nelson--Aalen estimators 
for the cumulative hazards $A_j(t)$ are given by
\begin{align*}
\widehat A_j(t)=\int_{0}^t \frac{\mathbf{1}\{Y_j>0\}}{Y_j}\,\mathrm{ d }N_j\quad (j=1,2).
\end{align*}

Before we subsequently define the weighted logrank test statistic 
we have to introduce the pooled Nelson--Aalen $\widehat A(t)=\int_0^t \mathbf{1}\{Y>0\}/Y \,\mathrm{ d }N$ $(t\geq 0)$ as well as the pooled Kaplan--Meier estimator $\widehat F$ given by
\begin{align*}
1-\widehat F(t)=\prod_{(j,i):X_{j,i}\leq t}\Bigl( 1- \frac{\delta_{j,i}}{Y(X_{j,i})} \Bigr)=\prod_{(j,i):X_{j,i}\leq t}\Bigl( 1- \frac{\Delta N(X_{j,i})}{Y(X_{j,i})} \Bigr)\quad (t\geq 0).
\end{align*}
Here, $\Delta N(t)=N(t)-N(t-)$ denotes the increment of the counting process $N$ at time $t$. Then the weighted logrank statistic \citep{FlemingETAL1987, ABGK} is defined as 
\begin{align*}
T_n(w)= \Bigr(\frac{n}{n_1n_2}\Bigl)^{1/2} \int_{0}^\infty w(\widehat F(t-)) \frac{Y_1(t)Y_2(t)}{Y(t)}\Bigl[ \,\mathrm{ d }\widehat A_1(t) -\,\mathrm{ d }\widehat A_2(t)\Bigr],
\end{align*}
where the weight function $ w\in \mathcal W$ corresponds to the main alternative of interest and is taken from the space 
\begin{align*}
\mathcal W =\{w:[0,1]\to[0, \infty)\text{ continuous and of bounded variation};\,w(x)> 0 \text{ for some }x\in[0,1]\}
\end{align*}  
in our one-sided case of interest. 
 \cite{Gill1980} proved asymptotic normality $T_n(w) \stackrel{d}{\to} N(0,\sigma^2)$ under the null hypothesis $H_=$ for each $w\in \mathcal W$. He suggested to estimate the limiting variance $\sigma^2$ by
\begin{align}\label{eqn:sigma_Gill}
\widehat\sigma_n^2 (w) = \frac{n}{n_1n_2} \int_{0}^{\infty} w(\widehat F(t-))^2 \frac{Y_1(t)Y_2(t)}{Y(t)} \,\mathrm{ d }\widehat A(t).
\end{align}
However, tests based on the one-sided test statistic $T_{n,stud}(w):=T_n(w)/\widehat\sigma_n (w)\ \mathbf{1}\{T_n(w)>0\}$ cannot be (asymptotically) optimal for the whole alternative $\kgeq$. Roughly speaking, 
the weighted logrank test $\phi_{n,\alpha}(w)=\mathbf{1}\{ T_{n,stud}(w) > z_{1-\alpha}\}$ 
leads to an 'optimal' decision for local alternatives of the form \citep{Gill1980}
\begin{align}\label{eqn:Kw}
	K_w: \left\{ {\lambda_1}/{\lambda_2}= 1+ \vartheta w(F_2)/\sqrt{n}: \vartheta> 0 
	\right\},
\end{align} 
see Section~\ref{sec:uncond_test} for its mathematically explicit form. Here, $z_{1-\alpha}$ denotes the $(1-\alpha)$-quantile and $\alpha\in(0,1/2)$ as we assume throughout. 
Restricting to $\vartheta>0$ in \eqref{eqn:Kw} ensures that $K_w\subset K_\geq$, i.e. stochastic ordering $F_1\geq F_2$, holds. As \eqref{eqn:Kw} is rather abstract we recall its meaning for the (classical) logrank test: 
As it is based on the constant weight function $w\equiv 1$ we thus obtain its optimality for proportional hazard alternatives of the form ${\lambda_1}/{\lambda_2}\equiv 1+ \vartheta/\sqrt{n}$. 
It is clear from \eqref{eqn:Kw} that  $K_w\subsetneq K_\geq$, i.e. the wrong choice of weight function $w$ may lead to a substantial loss in power. 
To enlarge $K_w$, \cite{BrendelETAL2014} adopted the projection idea of \cite{BehnenNeuhaus1989} to the present survival set-up, resulting in considerably larger optimality regions and broader power functions. We shortly explain their idea in the next Subsection and enhance it by means of a multiplier bootstrap technique and clever choice of weight functions.

\section{Asymptotics of Projection-type Tests}\label{sec:uncond_test}
For our investigations we throughout assume the usual asymptotic sample size condition 
\begin{enumerate}
 \item[(A1)] $n_1/n\to\kappa\in(0,1)$ as $\min(n_1,n_2)\to \infty$.
\end{enumerate}
Moreover, let  $\tau_j=\sup\{u\geq 0: [1-F_j(u)][1-G_j(u)]>0\}$ $(j=1,2)$
denote the observation times' upper limit within group $j$. 
To avoid the trivial case of purely censored data we additionally assume that $\max\{F_1(\tau),F_2(\tau)\}>0$  if  $\tau = \min\{\tau_1,\tau_2\}<\infty$.

\subsection{Optimal combination of weighted logrank statistics}

Let $w_1,\ldots,w_m \in \mathcal{W}$ be fixed weights corresponding to alternatives of interest (several examples are discussed in Sections~\ref{sec:simus}--\ref{sec:realdata} below). 
\cite{BrendelETAL2014} then proposed a test with broader power function on an enlarged alternative. Roughly speaking, their procedure is asymptotically optimal for $K_{w_1}, \dots, K_{w_m}$ as well as all local alternatives $K_w$ with arbitrary mixture of weights $w= \sum_{i=1}^{m}\beta_iw_i>0$ with $\beta_1,\ldots,\beta_m\geq0$. 
Denote their union as $K_m$. The choice of weights $w_1,\ldots,w_m\in \mathcal W$ together with the restriction to non-negative $\beta$'s ensures that $K_m\subset \kgeq$. 
It is, however, worth to note that the asymptotic results below do not require the constraint to weights from $\mathcal W$ but are even correct for more general  
continuous $w_\ell: [0,1]\to\R, w\neq 0$ of bounded variation, $\ell=1,\dots,m$. 

Anyhow, to achieve the outlined broader power behavior 
\cite{BrendelETAL2014} proceed as follows: 
Starting with the more simple two-sided test they pool all weighted logrank statistics in an $m$-dimensional vector ${\bf T}_n=[T_n(w_1),\ldots,T_n(w_m)]^T$ and use the quadratic form 
$Q_n := {\bf T}_n^T\widehat{\bf \Sigma}^-{\bf T}_n$ to detect two-sided deviations from the null hypothesis $H_=$.
Here, ${\bf B}^-$ denotes the Moore--Penrose inverse of a matrix ${\bf B}$ and $\widehat {\bf\Sigma}_n$ is the multivariate extension of \eqref{eqn:sigma_Gill} given by its entries 
\begin{align*}
	(\widehat {\bf\Sigma}_n)_{r,s}=\frac{n}{n_1n_2} \int_0^\infty w_s(\widehat F(t-))w_r(\widehat F(t-)) \frac{Y_1(t)Y_2(t)}{Y(t)} \,\mathrm{ d }\widehat A(t)\quad (r,s=1,\ldots,m).
\end{align*}
They then show that $Q_n$ can be expressed as an empirical projection of the hazard ratio difference
$\widehat{\Delta} = \mathrm d \widehat A_1/\mathrm d \widehat A-\mathrm d \widehat A_2/\mathrm d \widehat A$ 
onto the space $V=\{\sum_{i=1}^{m}\beta_iw_i(\widehat F):\beta_1,\ldots,\beta_m\in\R\}$ spanned by the pre-chosen weights $w_\ell$. In the one-sided case it is thus natural to restrict the projection 
of $\widehat{\Delta}$ to the positive cone $V_{\geq}=\{\sum_{i=1}^{m}\beta_iw_i(\widehat F):\beta_1,\ldots,\beta_m\geq 0\}$. They then prove that this eventually leads to the following maximum statistic in several quadratic forms
\begin{align}\label{eqn:def_S}
	S_n=\max\{0, {\bf T}_{n,J}^T\widehat{\bf \Sigma}_{n,J}^- {\bf T}_{n,J}: \emptyset\neq J \subset \{1,\ldots,m\};\,\widehat{\bf\Sigma}_J^-{\bf T}_{n,J}\geq 0\}.
\end{align}
Here, we used the notation ${\bf Z}_J=(Z_{j})_{j\in J}$ for a vector ${\bf Z}\in\R^m$ while ${\bf B}_J^-$ denotes the Moore-Penrose inverse of ${\bf B}_J=({\bf B}_{r,s})_{r,s\in J}$ for a matrix ${\bf B}\in \R^{m\times m}$ with indices taken from $J \subset \{1,\ldots,m\}$.

In the case $m=1$ of only one weight \eqref{eqn:def_S} equals the squared one-sided test statistic $S_n= T_n(w_1)^2/\widehat\sigma^2(w_1)\mathbf{1}\{ T_n(w_1)>0\}$ given below Equation \eqref{eqn:sigma_Gill}. Due to the rather complex structure of $S_n$ for general $m$, its asymptotic limit distribution is rather complicated and non-pivotal. 
 \cite{BrendelETAL2014} therefore propose a studentized permutation test in $S_n$ and prove its asymptotic optimality for testing $H_=$ against $K_m$. However, as outlined in the introduction, their procedure is computationally too exhaustive. A first step, to retort this issue is rather simple: We only treat linearly independent weights as this considerably simplifies the calculation of the involved Moore-Penrose inverse \citep{DitzhausFriedrich2018} and the limit distribution of $S_n$. As typical weights are given by polynomials (cf. Brendel {\it et al.}, 2014 and Sections \ref{sec:simus}-- \ref{sec:realdata} below), their linear independence is automatically given. Thus, the subsequent assumption, which we assume throughout, is no actual restriction from a practical point of view.
\begin{enumerate}
\item[(A2)] Suppose for all $\varepsilon\in(0,1)$ that $w_1,\ldots,w_m$ are linearly independent on $[0,\varepsilon]$, i.e., $\sum_{i=1}^m \beta_iw_i(x) = 0$ for all $x\in[0,\varepsilon]$ implies $\beta_1=\cdots=\beta_m=0$.
\end{enumerate}
Under this framework we obtain the following result about the asymptotic null distribution of $S_n$, which involves inverse matrices instead of general Moore-Penrose-inverses thanks to Assumption~A2.
\begin{theorem}[Convergence under the null]\label{theo:null_uncon}
	Suppose that the null $\ho$ is true and set $F_0=F_1=F_2$. 
	(a) The matrix $\widehat{{\bf \Sigma}}$ converges in probability to a non-singular matrix ${\bf \Sigma} =({\bf \Sigma}_{r,s})_{r,s}$ with entries ${\bf \Sigma}_{r,s}=\int (w_rw_s)(F_0)\psi\,\mathrm{ d }F_0$, where $\psi=[(1-G_1)(1-G_2)]/[\kappa(1-G_1)+(1-\kappa)(1-G_2)]$.\\
	(b) The test statistic $S_n$ in \eqref{eqn:def_S} converges in distribution to 
	\begin{align*}
		S=\max\{0,{\bf T}^T_J{\bf \Sigma}_J^{-1}{\bf T}_J: \emptyset\neq J \subset \{1,\ldots,m\};\,{\bf \Sigma}_J^{-1}{\bf T}_J\geq 0\},
	\end{align*}
	where ${\bf T}\sim N({\bf 0},{\bf \Sigma})$. Moreover, the distribution function $t\mapsto \P(S\leq t)$ of $S$ is continuous on $(0,\infty)$ with $\P(S=0)\leq 1/2$.
\end{theorem}
From Theorem \ref{theo:null_uncon} we obtain by  $\phi_{n,\alpha}=\mathbf{1}\{S_n>q_{m,1-\alpha}\}$ 
an asymptotically exact test, i.e., $\E(\phi_{n,\alpha})\to \alpha$, where $q_{m,1-\alpha}$ denotes the $(1-\alpha)$-quantile of the distribution of $S$. Due to Assumption~A2 the exact distribution of $S$ can be derived  as in the proof of Theorem 3.2.7 in \cite{BehnenNeuhaus1989}. 
However, it is rather cumbersome and depends on several unknown parameters and functions. To overcome this asymptotic non-pivotality of $S_n$ we later suggest an asymptotic correct wild bootstrap approach in Section \ref{sec:wild}. Before this, we state further asymptotic properties of the test $\phi_{n,\alpha}$ which will later carry over to its wild bootstrap version. 
To this end, we first note that in case of $m=1$ and $w_1=w$ our test is equivalent to the weighted logrank test $\phi_{n,\alpha}(w)=\mathbf{1}\{ T_{n,stud}(w) > z_{1-\alpha}\}$ based on $w$. 
With this in mind we discuss the consistency of our test $\phi_{n,\alpha}$ for general $m$ and explain how it combines the strength of these classical singly-weighted tests.
\begin{theorem}[Consistency]\label{theo:cons_uncon}
	Let $K\subset K_{\geq}$ be any fixed alternative. Then the test 
	$\phi_{n,\alpha}$ is consistent for testing $\ho$ versus $K$, i.e.,  $\E_K(\phi_{n,\alpha})\to 1$, whenever a singly-weighted logrank test $\phi_{n,\alpha}(w_i)$ is consistent for some $i=1,\ldots,m$.
\end{theorem}
For the singly-weighted logrank tests it was already shown \cite[Theorem 7.3.1]{FlemingHarrington} that $\phi_{n,\alpha}(w)$ with strictly positive weight $w>0$ is consistent for fixed ordered alternatives $F_1\geq F_2$ with $F_1(t)>F_2(t)$ for some $t<\tau$. 
Hence, we recommend to choose at least one strictly positive weight $w_i$, where we prefer to include the classical logrank weight $w\equiv 1$. 

To state the asymptotic optimality with respect to $K_m$ we now explicitly specify 
the corresponding local alternatives. To this end, let $F_0$ be a baseline distribution function with corresponding cumulative hazard function $A_0$ and consider the local alternatives
\begin{align}\label{eqn:local_alternative}
A_{j,n}(t)=\int_0^t 1+ c_{j,n} w(F_0) \,\mathrm{ d } A_0\quad (t\geq 0),\quad c_{j,n}= \frac{(-1)^{j+1}}{n_j}\Bigl( \frac{n_1n_2}{n} \Bigr)^{1/2}
\end{align}
for some $w\in \mathcal W$. Due to $w\geq 0$ we obtain a local alternative with $A_{1,n}\geq A_{2,n}$. 
Note, that different to the rather lax description of $K_m$ we here assume perturbations of $\lambda_j, j=1,2,$ in opposite directions which is needed to prove
\begin{theorem}[Asymptotics under local alternatives]\label{theo:local_altern}
	Define $\psi$ and ${\bf \Sigma}$ as in Theorem \ref{theo:null_uncon}, where now $F_0$ is the baseline distribution function. Then $S_n$ converges in distribution under the alternative \eqref{eqn:local_alternative} to a random variable $S$ as in \eqref{eqn:def_S}, where now
	$T\sim N({\bf a},{\bf \Sigma})$ is multivariate normal distributed with expectation vector ${\bf a}=( \int w(F_0)  w_i( F_0) \psi \,\mathrm{ d }F_0)_{i=1,\ldots,m}$. 
\end{theorem}

This result can be used to prove that $\phi_{n,\alpha}$ is admissible for all local alternatives of the form \eqref{eqn:local_alternative} with $w \in \mathcal W_m=\{w=\sum_{i=1}^m\beta_iw_i:\;\beta_1,\ldots,\beta_m \geq 0, w>0\}$, the cone spanned by the pre-chosen weights. 
In other words, there is no test achieving higher asymptotic power for all local alternatives given by weights $w\in\mathcal W_m$ simultaneously. 
To state this, denote by $Q_{n,{\pmb \beta}}$, ${\pmb \beta}=(\beta_1,\ldots,\beta_m)\in[0,\infty)^m,$ the distribution of the data $(X_{1,1},\delta_{1,1},\ldots,X_{2,n_2},\delta_{2,n_2})$ under the local alternative \eqref{eqn:local_alternative} with $w=\sum_{i=1}^m\beta_iw_i$. Furthermore, let $\E_{n,{\pmb \beta}}$ be the expectation under $Q_{n,{\pmb \beta}}$. In particular, $Q_{n,{\bf 0}}$ and $\E_{n,{\bf 0}}$ represent the situation under the null hypothesis $\ho$.

\begin{theorem}[Asymptotic admissible]\label{theo:local_admiss}
	There exists no test sequence $\varphi_n$ of asymptotic size $\alpha$, i.e. with $\limsup_{n\to\infty}\E_{n,{\bf 0}}(\varphi_n)\leq \alpha$, such that $\liminf_{n\to\infty}[\E_{n,{\pmb \beta}}(\varphi_n)-E_{n,{\pmb \beta}}(\phi_{n,\alpha})]$ is non-negative for all ${\bf 0}\neq {\pmb \beta} \geq 0$ and positive for some $w\in\mathcal W_m$.
\end{theorem}

In this sense $\phi_{n,\alpha}$ is asymptotically optimal for the alternative spanned by the pre-chosen weights $w_\ell$. How to calculate critical values efficiently, is discussed below.

\subsection{Wild Bootstrap}\label{sec:wild}
For a practical application of $\phi_{n,\alpha}$ critical values have to be calculated. As $q_{m,1-\alpha}$ depends on unknown quantities we propose a multiplier wild bootstrap technique \citep{Wu1986} which has previously been applied in other survival designs \citep{Lin1997, BeyersmannETAL2013, dobler2014bootstrapping, bluhmki2018wildNA}. To this end, let $G_{1,1},\ldots,G_{1,n_1}$, $G_{2,1},\ldots,G_{2,n_2}$ be $n$ independent and identical distributed multipliers with $\E(G_{j,i})=0$, $\var(G_{j,i})=1$ and finite fourth moment $E(G_{j,i}^4)<\infty$. The $G_{j,i}$ are also called wild bootstrap weights and are supposed to be independent of the data. As in \cite{bluhmki2018wildNA} the first idea is to randomly weight the subject-specific counting processes $N_{j,i}$ within the Nelson--Aalen estimator $\widehat A_j$ leading to its wild bootstrap version 
\begin{align}\label{eqn:def_AG}
	\widehat A_j^G(t)= \int_0^t \frac{\mathbf{1}\{Y_j>0\}}{Y_j}\,\mathrm{ d }\Bigl( \sum_{i=1}^{n_j}G_{j,i}N_{j,i} \Bigr)=\sum_{i=1}^{n_j}G_{j,i}\int_0^t \frac{\mathbf{1}\{Y_j>0\}}{Y_j}\,\mathrm{ d }N_{j,i}\quad(t \geq 0), 
\end{align}
$j=1,2.$
Now, replacing $\widehat A_j$ by $\widehat A_j^G$ in the definitions of $T_n(w)$ and ${\bf T}_n$ we obtain their wild bootstrap counterparts $T_n^G(w)$ and ${\bf T}_n^G$, respectively. This could already be used to obtain a wild bootstrap version of $S_n$. However, it is generally recommended to additionally bootstrap the studentization \citep{hall1991two}, i.e. the empirical covariance matrix $\widehat {\bf\Sigma}_n$. Following \citet{dobler2014bootstrapping, dobler2017non} and \cite{dobler2016nonparametric} a multiplier bootstrap version thereof is given by 
\begin{align}\label{eq:cov_boot}
	(\widehat {\bf\Sigma}^G_n)_{r,s}=
	\sum_{j=1}^2\sum_{i=1}^{n_j}
	\frac{n}{n_1n_2} \int_0^\infty w_s(\widehat F(t-))w_r(\widehat F(t-)) \frac{Y_1(t)Y_2(t)}{Y^2(t)} 
	G_{j,i}^2 \,\mathrm{ d } N_{j,i}(t)\quad (r,s=1,\ldots,m).
\end{align}
Now, substituting ${\bf T}_n$ by ${\bf T}_n^G$ and 
$\widehat {\bf\Sigma}_n$ by $\widehat {\bf\Sigma}^G_n$ in the definition of $S_n$ leads to its wild bootstrap version $S_n^G$. 
We note that another possibility for bootstrapping $\widehat {\bf\Sigma}_n$ as, e.g., proposed by \cite{BeyersmannETAL2013}, is given by the empirical covariance matrix of ${\bf T}_n^G$. However, our simulation results preferred the usage of \eqref{eq:cov_boot}, cf. Section~\ref{sec:simus} below and the appendix for details.

Now, only the question remains which weights should be chosen. Due to the asymptotic normality of ${\bf T}_n$ it seems plausible to choose standard normal $G_{j,i}$ as proposed by \cite{Lin1997}. Regarding the discrete character of the involved counting processes, however, it might also be 
reasonable to utilize discrete distributed multipliers instead to better reflect their finite sample behavior. This has, e.g., been argued by \cite{BeyersmannETAL2013} and \cite{dobler2017non} who proposed the use of centered Poisson multipliers. Another possibility is given by simple random signs (so-called Rademacher weights) that are uniformly distributed on $\{-1,1\}$. In case of linear regression models the latter exhibited preferably finite sample properties \citep{Liu1988,davidson2008wild}. Anyhow, we formulate the asymptotic validity of this approach for arbitrary weights and thereafter compare the three above multiplier choices in simulations.

\begin{theorem}[Wild bootstrap under the null]\label{theo:wild_boot_null}
	Define $S$ as in Theorem \ref{theo:null_uncon}. Then the wild bootstrap version $S^G_n$ asymptotically mimics its null distribution, i.e. under $\ho$ we have convergence in probability
	\begin{align}\label{eqn:theo_wild_boo_state}
		\sup_{x\geq 0}\Bigl | \P( S_n^G\leq x \vert (X_{1,1},\delta_{1,1}),\ldots,(X_{2,n_2},\delta_{2,n_2}) ) - \P(S\leq x) \Bigr| \to 0.	
	\end{align}
\end{theorem}
Consequently, $\phi_{n,\alpha}^G= \mathbf{1}\{S_n^G>q_{n,1-\alpha}^G\}$ is an asymptotically exact level $\alpha$ test for $\ho$, where  $q_{n,\alpha}^G$ 
denotes the $\alpha$-quantile of $S_n^G$ given the data $(X_{j,i},\delta_{j,i})_{1\leq i \leq n_j,j=1,2}$. 
Different to the permutation approach of \cite{BrendelETAL2014} the proposed multiplier resampling technique directly recovers  the asymptotic null distribution of the test statistic. This is not only methodologically desirable but also allows to work without any further modifications of the resampling as, e.g., a time-consuming studentization technique as needed in \cite{BrendelETAL2014}.

In addition to the asymptotic behavior under the null, the power behavior of $\phi_{n,\alpha}$ under fixed and local alternatives carries over to its wild bootstrap version $\phi_{n,\alpha}^G$. 
\begin{theorem}[Wild bootstrap under fixed and local alternatives]\label{theo:wild_boot_alt}
	(a) Let $K\subset K_{\geq}$ be any fixed alternative. Then $\phi_{n,\alpha}^G$ is consistent for $K$ whenever $\phi_{n,\alpha}$ is consistent for $K$.\\
	(b) Under the local alternatives \eqref{eqn:local_alternative} $\phi_{n,\alpha}^G$ and $\phi_{n,\alpha}$ share the same asymptotic power behavior, i.e., $\E_{n,w}(|\phi_{n,\alpha}-\phi_{n,\alpha}^G|)\to 0$. In particular, $\phi_{n,\alpha}^G$ is asymptotically admissible in the sense of Theorem \ref{theo:local_admiss}.
\end{theorem}

\section{Simulations}\label{sec:simus}

To complement our theoretical large sample studies from above a simulation study was conducted in which we investigate our tests' small sample properties. 

\subsection{Type-I error}\label{sec:simus_null}

Since monotonic transformation of the observations do not affect the rank tests' outcome we can restrict ourselves to standard exponential distributed survival times. The censoring times are simulated by exponential distributions $\text{Exp}(\mu_j)$, where the scale parameter $\mu_j>0$ $(j=1,2)$ depends on the group only. In the first setting, the parameter $\mu_j$ is chosen such that it reflects an average censoring rate of $10\%$ in the first group and of $30\%$ in the second group. In addition, we also discuss two cases of equal censoring with rates of $15\%$ and $30\%$, respectively. Moreover, we consider four different scenarios $(n_1,n_2)=$ $(25,25)$, $(20,30)$, $(50,50)$, $(30,70)$ representing balanced and unbalanced settings with small to moderate sample sizes. The nominal level $\alpha$ is set to $5\%$. After playing around with several weights we consider three directions given by
\begin{align}\label{eqn:def_w}
	w_1(x) = 1, \quad w_2(x) = (1-x)^4, \quad w_3(x) = x^4\quad (x\in[0,1]).
\end{align}
Here, $w_1$ corresponds to the classical logrank test giving equal weight to all times. In contrast, $w_2$ and $w_3$ 
put more emphasize on early and late times, respectively. It is easy to check that the weights are linearly independent. Thus,  Assumption $A2$ is fulfilled. We also investigate the effect of different wild bootstap multipliers. In particular, we 
 choose from either Rademacher, normal or centered Poisson multipliers $G_{j,i}$. 
All computations are done using the computing environment \textsc{R} \citep{R}, version 3.5.0, based on $N_{sim}=5,000$ simulation runs and $N_{boot}=1,000$ bootstrap runs for each scenario.
The resulting empirical sizes are displayed in Table \ref{table:null}.  

\begin{table}[ht]
	\caption{Empirical sizes in $\%$ (nominal level $5\%$) for three different wild bootstrap versions of our test. The survival and censoring times were simulated by exponential distributions.  }\label{table:null}
	\centering
	\begin{tabular}{ccccc}
		\hline
		$(n_1,n_2)$ & censoring in \% & Normal & Poisson & Rademacher  \\ 
		\hline \hline
				   & (10,30) &  6.02 & 9.14 & 5.18 \\ 
		(20,30)    & (15,15) & 5.96 & 8.98 & 4.66 \\ 
		           & (30,30) & 5.86 & 8.66 & 4.24 \\  
		\hline \hline
		           & (10,30) & 5.54 & 9.04 & 5.14 \\ 
		 (25,25)   & (15,15) & 5.84 & 9.02 & 5.16 \\ 
		           & (30,30) & 5.86 & 9.06 & 5.46 \\
		 \hline\hline
		 	       & (10,30) & 5.38 & 8.28 & 4.72 \\ 
		 (30,70)   & (15,15) & 5.02 & 8.04 & 4.40 \\  
		           & (30,30) & 5.92 & 8.62 & 4.38 \\ 
		 \hline \hline
		 		   & (10,30) & 5.48 & 8.50 & 5.72 \\  
		 (50,50)   & (15,15) & 5.80 & 8.84 & 5.38 \\ 
		 	       & (30,30) & 5.26 & 8.24 & 5.10 \\ 
	\end{tabular}
\end{table}

It is apparent that the wild bootstrap test $\phi_{n,\alpha}^G$ based on centered Poisson multipliers tends to quite liberal conclusions with true type-$I$-error rates around $8 - 9 \%$. This is completely different to findings in competing-risks settings \citep{BeyersmannETAL2013}, where centered Poisson appeared to be preferable over standard normal ones. In particular, in our situation choosing normal or Rademacher multipliers considerably enhances the type-$I$-error control with rates between $4 - 6 \%$. Thereof, we recommend Rademacher multipliers due to a slight liberality of the test based on normal ones. Simulation results for empirical covariance estimators 
of ${\bf T}^G_n$ instead of $\widehat {\bf\Sigma}^G_n$ given in the appendix did not alter this general conclusion. However, the test based on Poisson multipliers then exhibits a better behavior with type-$I$-errors between $5-8\%$ while the two others tend to be more liberal in this case. 
We therefore restrict our further power investigations to Rademacher multipliers together with $\widehat {\bf\Sigma}^G_n =  \widehat {\bf\Sigma}_n$ (in this case) as covariance estimator in the bootstrap test.


\subsection{Power behavior against various alternatives}\label{sec:simu_alternative}
\begin{figure}[h!] 
	\begin{center}	
		\includegraphics[width=0.9\textwidth]{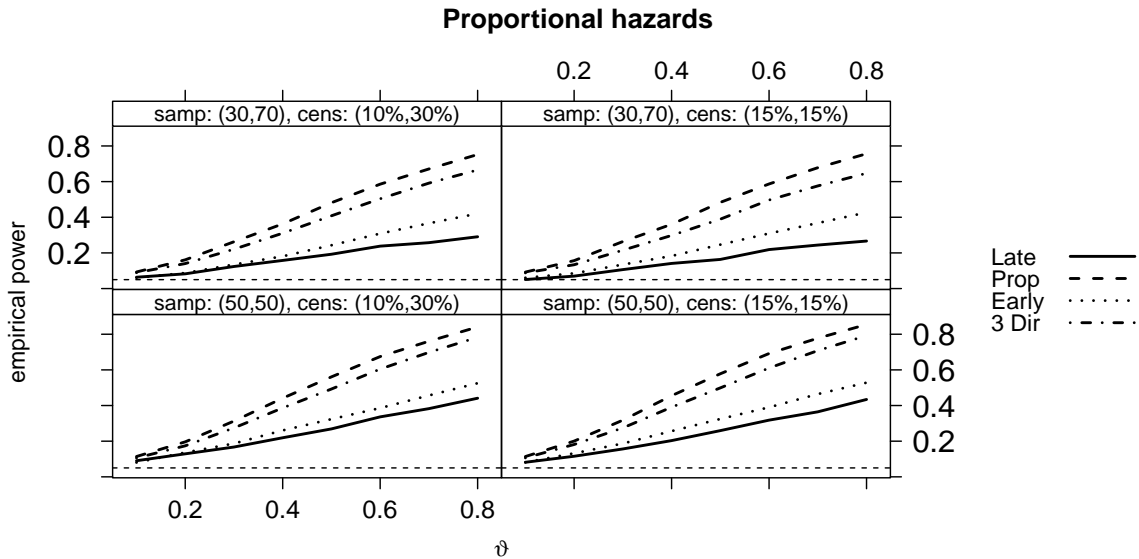}
	\end{center}
	\caption[h]{Power simulation results ($\alpha = 5\%$) under the proportional hazard alternative for the test based on all three directions (3 Dir) and the three singly-weighted tests based on one of these directions.} \label{fig:prop}
\end{figure}

\begin{figure}[h!] 
	\begin{center}	
		\includegraphics[width=0.9\textwidth]{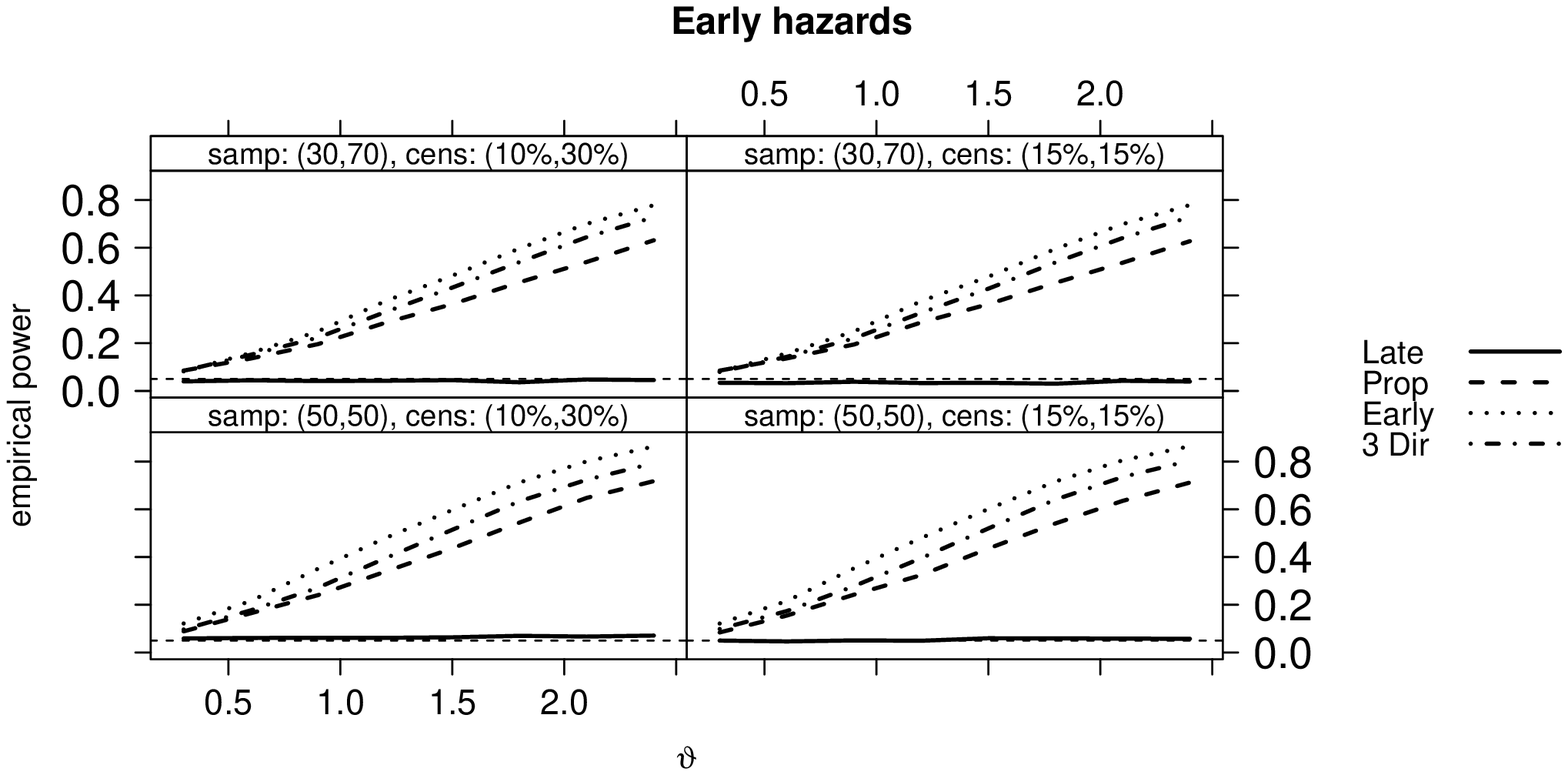}
	\end{center}
	\caption[h]{Power simulation results ($\alpha = 5\%$) under the early hazard alternative for the test based on all three directions (3 Dir) and the three singly-weighted tests based on one of these directions.} \label{fig:early}
\end{figure}

\begin{figure}[h!] 
	\begin{center}	
		\includegraphics[width=0.9\textwidth]{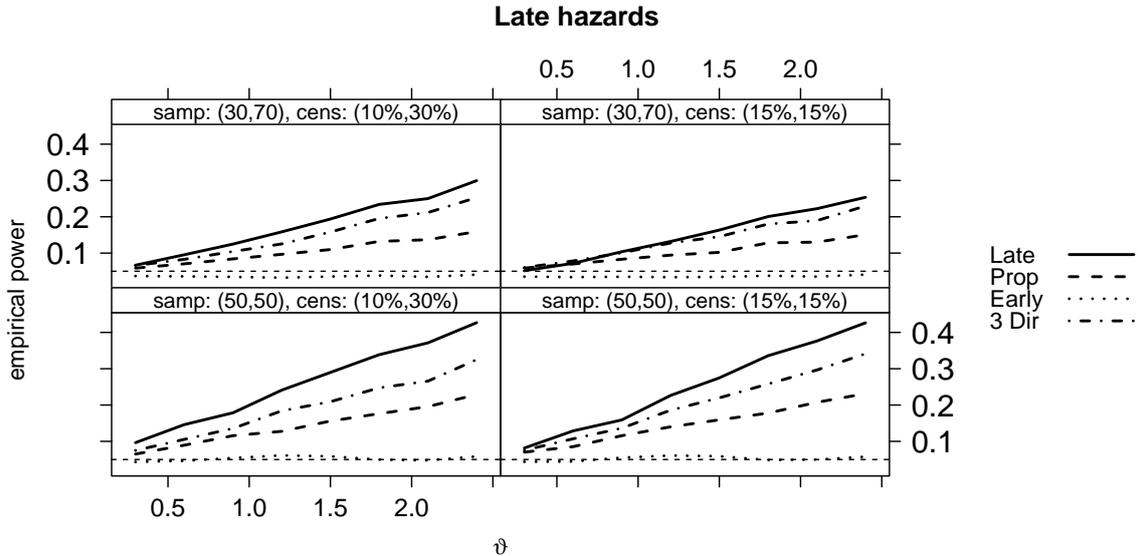}
	\end{center}
	\caption[h]{Power simulation results ($\alpha = 5\%$) under the late hazard alternative for the test based on all three directions (3 Dir) and the three singly-weighted tests based on one of these directions.} \label{fig:late}
\end{figure}

In this section we present simulation results for our test's power under various scenarios. Regarding the results of Section \ref{sec:simus_null} we consider only the Rademacher wild bootstrap version of our test based on $w_1,w_2,w_3$ from \eqref{eqn:def_w}. We compare our test with the Rademacher wild bootstrap versions of the singly-weighted tests based on only one of the three weights/directions, respectively. The simulation settings are the same as in Section \ref{sec:simus_null} except for the first group's distribution. Here, we replace the standard exponential distribution by a perturbation of it in the three hazard directions $w_1$ (proportional), $w_2$ (early) and $w_3$ (late). To be more specific, we consider the situation described in \eqref{eqn:Kw} for $w=w_i$: We choose $8$ equidistant values of $\vartheta$ from $0.3$ to $2.4$ for the early ($w_2$) and late ($w_3$) hazard alternatives, respectively, as well as from $0.1$ to $0.8$ for the proportional ($w_3$) hazard alternative. For ease of presentation we only consider the larger sample size settings $(n_1,n_2)=$ $(50,50)$, $(30,70)$ with censoring rates $(15\%,15\%)$ and $(10\%,30\%)$, respectively. The tests' empirical power curves are displayed in Figures \ref{fig:prop}--\ref{fig:late}.

In all scenarios the same phenomenon can be readily seen: The singly-weighted test $\phi_{n,\alpha}(w)$ based on the same weight $w$, which is used for the alternative, leads to the highest empirical power values in each scenario, e.g.
$\phi_{n,\alpha}(w_2)$ had largest power for early hazards. Due to the admissibility result stated in Theorems \ref{theo:local_admiss} and \ref{theo:wild_boot_alt} this is not surprising. The test $\phi_{n,\alpha}$ based on the three directions $w_1,w_2,w_3$ always leads to the second highest values. This confirms the primary intention and its strength: instead of having an optimal test for a single weight alternative we broaden the power and get a test with good power behavior under various alternatives. Having this goal in mind one may choose far more than three weights. But, as already pointed out by \cite{DitzhausFriedrich2018} for the two-sided situation, choosing too many weights may result in flatter power curves. In addition, this also increases the computational load of the procedure 
as the maximum statistic \eqref{eqn:def_S} requires the calculation of $2^m-1$ (sub-) quadratic forms.

\begin{figure}[h!] 
	\begin{center}	
		\includegraphics[width=0.9\textwidth]{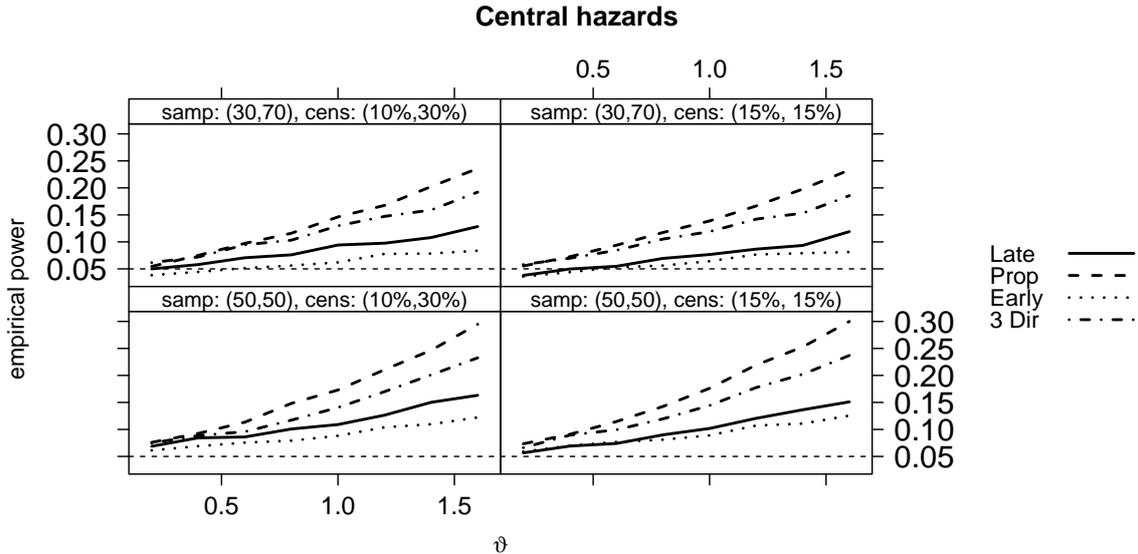}
	\end{center}
	\caption[h]{Power simulation results ($\alpha = 5\%$) under the central hazard alternative for the test based on all three directions (3 Dir) and the three singly-weighted tests based on one of these directions.} \label{fig:central}
\end{figure}

Finally, we also consider a scenario for which the alternative does not correspond to the above weights $w_1,w_2,w_3$. We pick  $w(x)=x(1-x)$ giving more weight to 'central' times, i.e., times around the median, and choose $8$ equidistant values of $\vartheta$ between $0.2$ and $1.6$. Here, one would a priori suspect a worse power behavior of the singly-weighted logrank tests which are sensitive to early and late hazard alternatives. In fact, the results displayed in Figure \ref{fig:central} are in line with our intention and the previous power results. 

Moreover, the singly-weighted classical logrank test $\phi_{n,\alpha}^G(w_1)$ has the highest power values in all scenarios and the novel test $\phi_{n,\alpha}^G$ again exhibits the second highest power; closely following the power curve of $\phi_{n,\alpha}^G(w_1)$.

\section{Real data example}\label{sec:realdata}
\begin{figure}[ht] 
	\begin{minipage}{0.5\textwidth}
		\includegraphics[width=\textwidth]{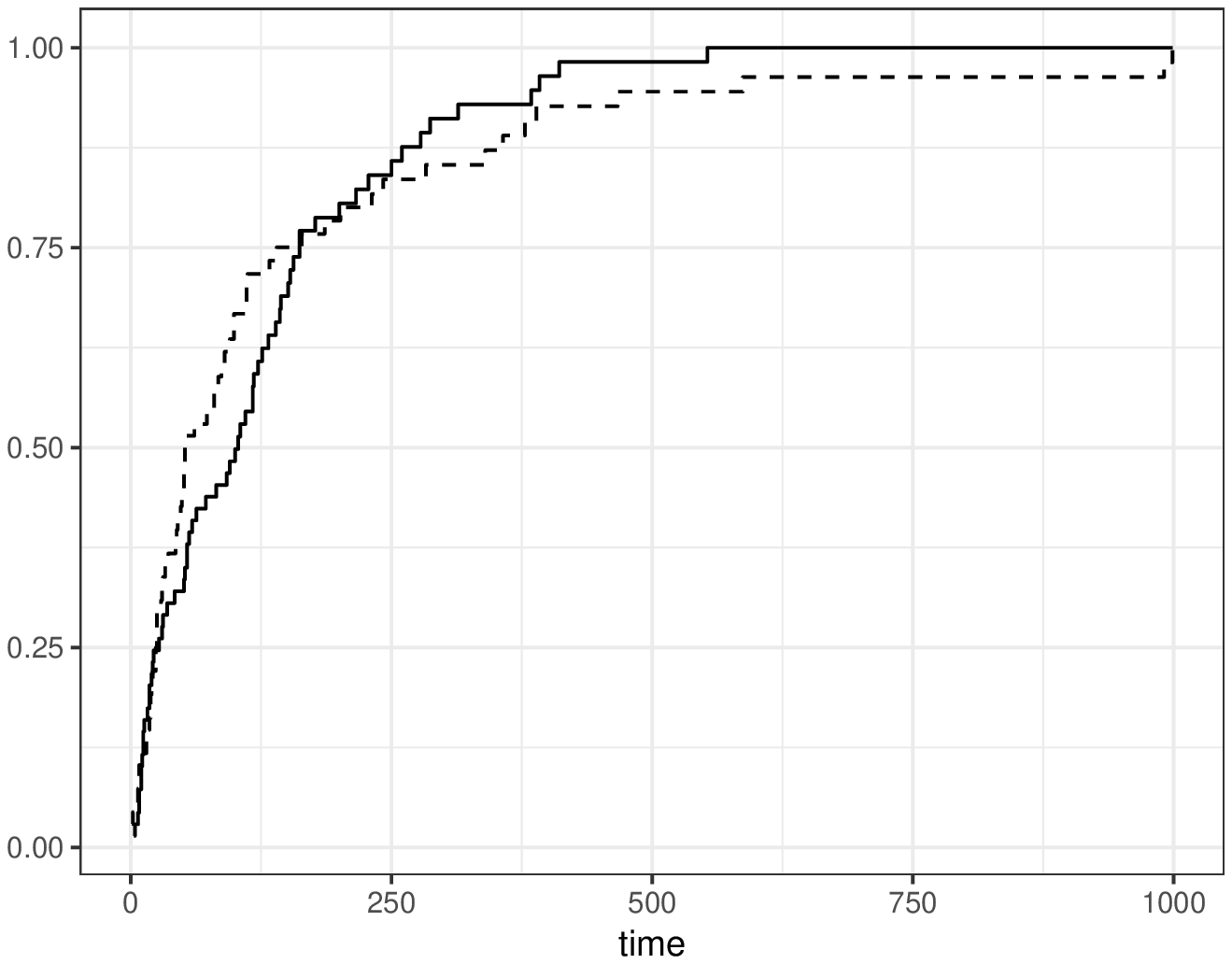}
	\end{minipage}
	\begin{minipage}{0.5\textwidth}
		\includegraphics[width=\textwidth]{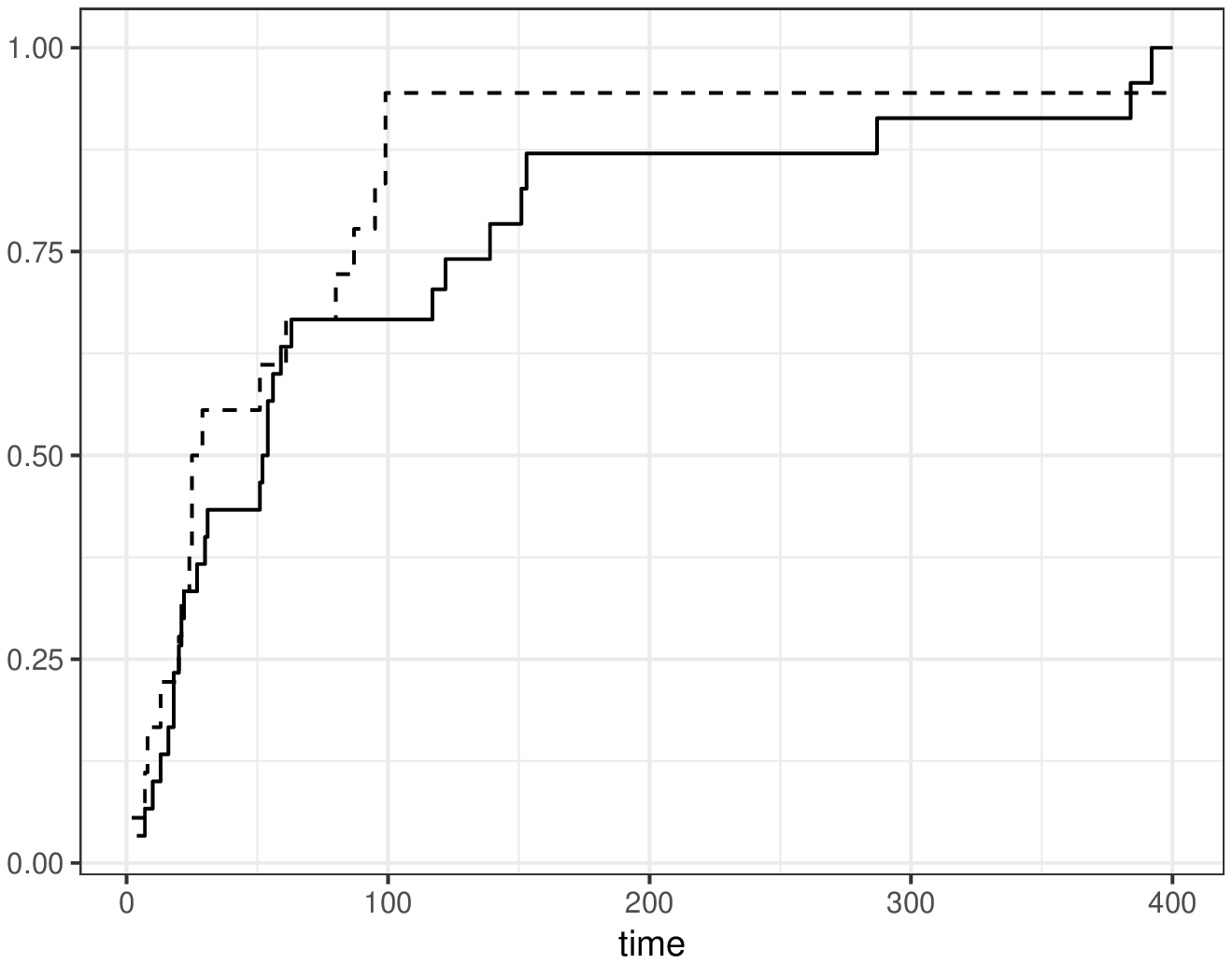}
	\end{minipage}
	\caption[h]{(Group-specific) Kaplan--Meier curves of the veteran data set considering all tumor types (left) or only the type 'small cell' (right). The curves corresponding to the standard chemotherapy are solid (---) and the ones of the experimental chemotherapy are dashed (- - -), respectively.} \label{fig:KMcurves}
\end{figure}

To illustrate the applicability of our procedure we re-analyze the Veteran's Administration lung cancer study from \cite{Prentice}. The data set is included in the \textsc{R}-package \textit{survival} via the command \textit{veteran}. It  consists of $137$ survival times of males with inoperable lung cancer, including $9$ censored observations corresponding to a censoring rate of nearly $7\%$. In the study two interventions, a standard and an experimental chemotherapy, were compared on four different tumor types. As  the tumor type had a statistically significant impact  on the survival time 
\citep{KalbfleischPrentice1980} we restrict our exemplary analysis to the tumor type \textit{small cell}. In this setting, group 1 corresponds to the experimental chemotherapy while the standard chemotherapy is represented by group 2. We test the alternative $K_\geq$ that the experimental chemotherapy has a negative effect on the patients' mortality compared to the standard treatment. Although it contradicts the findings of \cite{KalbfleischPrentice1980}, we also analyze the whole data set ignoring the different tumor types for didactic reasons; now testing superiority of the experimental group. The (group-specific) Kaplan-Meier curves are displayed in Figure \ref{fig:KMcurves} for both examples.

To check $K_\geq$ we apply our test based on the early, proportional and late weights $w_1,w_2,w_3$ as already done in the previous section. As recommended in Section \ref{sec:simus_null} we only consider the Rademacher wild bootstrap multipliers. To estimate the corresponding bootstrap quantile we use $N_{boot}=10^4$ iterations. In particular, these are exactly the default configurations of Marc's novel \textsc{R}-function \textit{mdir.onesided} implemented within the package {\bf mdir.logrank} available on GitHub. We briefly illustrate its application on the current study: After loading the data set the columns consisting of the survival times, the censoring and the group  status are renamed to \textit{time}, \textit{event} and \textit{group}, respectively. 
For example, in the veteran data set the coding $1$ for standard and $2$ for experimental chemotherapy was used (i.e., the group status is permuted when compared to ours). Then the command \textit{mdir.onesided(data, group1 = 2)} calculates the desired $p$-value of a superiority of the standard chemotherapy group (coded as group~1 in the \textsc{R} data set). Hence, \textit{mdir.onesided(data, group1 = 1)} would lead to a test for superiority of the experimental group. 
 For the user's convenience we also implemented a GUI, where the user can easily specify weights of the form $w(x)=x^r(1-x)^g$ with $r,g\in\N$ or even choose one of the other two wild bootstrap approaches. 
Thereby, it is automatically checked whether Assumption $A2$ is fulfilled for the chosen weights. In case of linearly dependent weights a linearly independent subclass is selected automatically. The GUI is called upon the command \textit{calculateGUI()} and a detailed step-by-step description of its applicability with corresponding screen-shots is given in the appendix. Finally, more general weights are also accessible but only via the pure command \textit{mdir.onesided}. 

In addition to our test, we calculate the $p$-values of the classical singly-weighted tests $\phi_{n,\alpha}(w_i)$ for $i=1,2,3$. The resulting $p$-values are displayed in Table \ref{table:real_data_pvalues}. The plot of the Kaplan--Meier curves for the first example suggest a superiority of the first group's distribution function, or in other words the alternative seem to be true. This first graphical impression can be supported statistically since our test's $p$-value is smaller than $5\%$. Contrary, the application of the classical logrank test $\phi_{n,\alpha}(w_1)$ or the test $\phi_{n,\alpha}(w_2)$ based on an early weight would not lead to a significant rejection. For the pooled data set (consisting of all tumor types), $\phi_{n,\alpha}(w_2)$ again rejects the null, whereas $\phi_{n,\alpha}(w_1)$ and $\phi_{n,\alpha}(w_2)$ do not reject. This is in line with Figure~\ref{fig:KMcurves} since the Kaplan--Meier estimator $\widehat F_1$ for the first group exceeds the one $\widehat F_2$ of the second group for late times but is smaller than $\widehat F_2$ for early times. In all, no overall superiority can be graphically proven, which is also the statistical result of our test. 
\begin{table}[ht]
	\caption{$p$-values of our test $\phi_{n,\alpha}^G$ based on early, proportional and late hazard weight as well as of the classical singly-weighted tests $\phi_{n,\alpha}(w_i), i=1,2,3,$ for the veteran data set respecting all tumor types or only the type small cell, respectively.}\label{table:real_data_pvalues}
	\centering
	\begin{tabular}{rrrrr}
		\hline
		data set & 3 weights & prop & early & late \\ 
		\hline \hline
		small cell tumor & 0.043 & 0.066   &  0.279  &  0.003 \\
		all tumor types  & 0.086 & 0.533   &  0.703  &  0.028      
	\end{tabular}
\end{table}

\section{Discussion and Outlook}\label{sec:discussion}
We investigated the problem of testing for superiority in an unpaired two-sample survival set-up, allowing for heterogeneous censoring structures. 
To enlarge the power functions of classical weighted logrank tests, we followed the idea of \cite{BrendelETAL2014} and enhanced it in two directions. 
They proposed a projection-type permutation test that is asymptotically optimal for a fixed number of pre-determined alternatives such as early, proportional or late hazard differences. However, their procedure is computationally exhaustive; particularly due to a complicated permutation technqiue that accounts for the test statistics' complex 
limit distribution and also needs two time-consuming Monte-Carlo loops for calculating $p$-values. To this end, we first simplified their statistic by only allowing the combination of alternatives corresponding to linearly independent weight functions. We explain that this is no loss at all from a practical point of view but instead leads to a more simple limit distribution. This allows for a more simple and direct approximation by means of wild bootstrap multiplier resampling. The so gained wild bootstrap test is then shown to be asymptotically correct under the null, consistent and also optimal for the chosen alternatives of interest. To give recommendations for the choice of multipliers simulations for several weights corresponding to alternatives of interest were run. Instead of the usual standard normal we found a preference for Rademacher multipliers. Moreover, to illustrate the effect of the chosen weights corresponding to alternatives of interest we conducted several power simulation resulting in the simplified recommendation to include at least weights for early, proportional and late hazard differences. 

The new methods are implemented in the \textsc{R}-package {\bf mdir.logrank} available on GitHub, where Rademacher multipliers and a combination of 
weights for early, proportional and late hazard alternatives are set as the default choice. However, different choices are also possible and  implemented. An analysis of a two-armed lung study explains its aplicatoin and additionally illustrates the effect of the different choices of weights.

The potential transfer of the current methodology to more complex designs, e.g. allowing for more than two groups or even more than the two states may be part of future research. 

\nocite{dobler2017bootstrap}

\section*{Acknowledgement}
This work was supported by the \textit{Deutsche Forschungsgemeinschaft}.


\bibliographystyle{plainnat}
\bibliography{sample}

\appendix

\section{ The Proofs}

For a slightly abbreviated formulation of the proof steps we set $\kappa_1=\kappa$ and $\kappa_2=1-\kappa$.

\subsection{Proof of Theorem \ref{theo:null_uncon}}

In the proof of their Theorem 1 \cite{DitzhausFriedrich2018} showed that the covariance matrix ${\bf \Sigma}$ is non-singular under Assumption $A2$. Moreover, \cite{BrendelETAL2014} already verified that $\widehat{{\bf \Sigma}}$ converges in probability to ${\bf \Sigma}$ and that ${\bf T}_n$ tends in distribution to ${\bf T}\sim N({\bf 0},{\bf\Sigma})$. Since ${\bf\Sigma}$ and so ${\bf\Sigma}_j$ is invertible for any subset $\emptyset \neq J\subset\{1,\ldots,m\}$ we can deduce immediately the convergence of the corresponding Moore--Penrose inverses $\widehat {\bf\Sigma}_J^-$ to ${\bf\Sigma}_J^-={\bf\Sigma}_J^{-1}$ as there finally is no rank jump in probability. Consequently, the convergence of $S_n$ to $S$ in distribution follows from the continuous mapping theorem.  The statements about the distribution function $t\mapsto\P(Z\leq t)$ and the probability $\P(Z=0)$ were already shown by \cite{BrendelETAL2014}, see Lemma 9.3 in their supplement.\qed

\subsection{Proof of Theorem \ref{theo:cons_uncon}}
Fix $i\in\{1,\ldots,m\}$. Clearly, $S_n\geq T_n(w_i)^2/\widehat{\sigma}^2(w_i)\mathbf{1}\{T_n(w_i)>0\}$. Hence, $\phi_{n,\alpha}$ is consistent for a fixed alternative $K$ whenever $\phi_{n,\alpha}(w_i)$ is consistent for $K$. This implication was already observed by \cite{BrendelETAL2014}, see their Theorem 2.\qed

\subsection{Proof of Theorem \ref{theo:local_altern}}
Combining Theorem 7.4.1 of \cite{FlemingHarrington} and the Cram\'{e}r--Wold device we can verify the distributional convergence of ${\bf T}_n$ to ${\bf T}$. This was already done by \cite{BrendelETAL2014}, see their  Theorem 9.1 in the supplement. The final statement of Theorem \ref{theo:local_altern} follows from the continuous mapping theorem. \qed

\subsection{Proof of Theorem \ref{theo:local_admiss}}
As already indicated in the main paper our test statistic can be interpreted as a certain projection, see Theorem 1 of \cite{BrendelETAL2014}.  \cite{BehnenNeuhaus1989} studied this type of projection statistics and verified the following representation of $S_n$, see their Equations (3.2.1) and (6.3.3):
\begin{align}\label{eqn:SN_representation}
S_n= \sup\{2 {\pmb\beta}^T{\bf T}_n - {\pmb\beta}^T\widehat {\bf\Sigma} {\pmb \beta}: {\pmb\beta}\in[0,\infty)^m\}.
\end{align}
Let $\mathcal M_m\subset \R^{m\times m}$ be the set of all positive definite matrices. Then $
\R^m \times \mathcal M_m \ni ({\bf x},{\bf \Gamma})\mapsto f({\bf x},{\bf\Gamma})= \sup\{ 2 {\pmb\beta}^T {\bf x} - {\pmb\beta}^T {\bf\Gamma} {\pmb\beta}: {\pmb\beta}\in [0,\infty)^m\}$ is jointly continuous in both arguments, see Lemma 7.5.7 of \cite{BehnenNeuhaus1989}. Moreover, ${\bf x}\mapsto f({\bf x},{\bf\Gamma})$ is convex for each fixed matrix ${\bf\Gamma}\in \mathcal M_m$. Both properties are needed in the following derivations. 

Let $Q_{{\pmb\beta}}=N({\bf\Sigma{\pmb\beta}},{\bf\Sigma})$ $({\pmb\beta}\in \R^m)$. Then it is easy to see that $\frac{\mathrm dQ_{n,{\pmb\beta}}}{\mathrm dQ_{n,{\bf 0}}}$ converges in distribution under $Q_{n,{\bf 0}}$ to $\frac{\mathrm dQ_{{\pmb\beta}}}{\mathrm dQ_{{\bf 0}}}({\bf Z})$ with ${\bf Z}\sim Q_{\bf 0}$ \citep[Proof of Theorem 4]{DitzhausFriedrich2018}. 
Using the notation of statistical experiments \citep[Sections 60 and 80]{Strasser1985}, this means that the statistical experiment $\{Q_{n,{\pmb\beta}}:{\pmb\beta}\in\R^m\}$ converges weakly to the Gaussian experiment $\{Q_{{\pmb\beta}}:{\pmb\beta}\in\R^m\}$. In other words, $\{Q_{n,{\pmb\beta}}:{\pmb\beta}\in\R^m\}$  fulfills Le Cam's asymptotic normality. We can now deduce from Le Cam's first lemma \citep[ Theorem 61.3]{Strasser1985} that $Q_{n,{\pmb\beta}}$ and $Q_{n,{\bf 0}}$ are mutual contiguous for all ${\pmb \beta}\in\R^m$. In particular, any convergence holds in $Q_{n,{\pmb \beta}}$-probability if and only if it does so in $Q_{n,{\bf 0}}$-probability. Combining this and the convergence of $\widehat {\bf \Sigma}$ under the null (Theorem \ref{theo:null_uncon}) yields that $\widehat {\bf \Sigma}$ converges in $Q_{n,{\pmb\beta}}$-probability to $\bf \Sigma$ for every ${\pmb\beta}\in \R^m$. Hence, we can deduce from Theorem \ref{theo:local_altern}, the continuity of $f$ and Lebesgue's dominated convergence theorem that
\begin{align*}
E_{n,{\pmb\beta}}(\phi_{n,\alpha}) = E_{n,{\pmb\beta}}(\mathbf{1}\{f({\bf T}_n,\widehat {\bf \Sigma})>q_{m,1-\alpha}\})\to \int \mathbf{1}\{f(x,{\bf \Sigma})>q_{m,1-\alpha}\} \,\mathrm{ d }Q_{{\pmb\beta}}(x).
\end{align*}
Define $\phi_{\alpha}^*= \mathbf{1}\{f(x, {\bf \Sigma})>q_{m,1-\alpha}\}$. The acceptance region $A=\{x\in\R^m: f(x,{\bf \Sigma})\leq q_{m,1-\alpha}\}$ of $\phi_{\alpha}^*$ is convex. Hence, we can conclude from Stein's theorem, see e.g. Theorem 5.6.5 of \cite{Anderson2003}, that $\phi_{\alpha}^*$ is admissible in the limiting model $\{Q_{{\pmb\beta}}:{\pmb\beta}\in[0,\infty)^m\}$ for testing $H:\{{\pmb\beta} = 0\}$ versus $K:\{{\pmb\beta}\geq 0,\,{\pmb\beta}\neq 0\}$. Contrary to our claim let us now assume, that there is a test sequence $\varphi_n$ $(n\in\N)$ of asymptotic size $\alpha$ such that $\limsup_{n\to\infty}E_{n,{\pmb\beta}}(\varphi_n-\phi_{n,\alpha})$ is nonnegative for all $0\neq {\pmb\beta}\geq 0$ and positive for at least one of these ${\pmb\beta}$'s. By a theorem of Le Cam \citep[Theorem 62.3]{Strasser1985} there is a test $\varphi$ of the limiting model $\{Q_{{\pmb\beta}}:{\pmb\beta}\in[0,,\infty)^m\}$ such that $\E_{n,{\pmb\beta}}(\varphi_{n})\to \int \varphi \,\mathrm{ d }Q_{{\pmb\beta}}$ for all ${\pmb\beta}\geq 0$ along an appropriate subsequence. But this implies that $\varphi$ is of size $\alpha$ for the limiting null $H:{\pmb\beta}=0$ and $\int \varphi -\phi_{\alpha}^* \,\mathrm{ d }Q_{{\pmb\beta}}$ is nonnegative for all $0\neq{\pmb\beta}\geq 0$ as well as positive for at least one of these ${\pmb\beta}$'s. Clearly, this contradicts the admissibility of $\phi_{\alpha}^*$ and, finally, proves the (asymptotic) admissibility of our test.\qed

\subsection{Proof of Theorems \ref{theo:wild_boot_null} and \ref{theo:wild_boot_alt}}
As explained below it is sufficient for all statements of Theorems \ref{theo:wild_boot_null} and \ref{theo:wild_boot_alt} to proof that  \eqref{eqn:theo_wild_boo_state} holds for some real-valued random variable $S_0$, say, under $\ho$ as well as under fixed alternatives $K$, respectively, where the distribution of $S_0$ depends on the underlying distributions $F_1,F_2,G_1,G_2$ and equals the one of $S$ under $\ho$. From the latter the asymptotic exactness under the null follows immediately. Moreover, we obtain $\E_{n,{\bf 0}}(|\phi_{n,\alpha}-\phi_{n,\alpha}^G|)\to 0$ from Lemma 1 of \cite{JanssenPauls2003}. Due to the mutually contiguity of $Q_{n,{\bf 0}}$ and $Q_{n,{\pmb\beta}}$ explained in the proof of Theorem  \ref{theo:local_admiss} this convergence carries over to local alternatives, i.e., we obtain $\E_{n,{\pmb\beta}}(|\phi_{n,\alpha}-\phi_{n,\alpha}^G|)\to 0$ for all ${\pmb\beta}\in[0,\infty)^m$. If $\phi_{n,\alpha}$ is consistent for a fixed alternative $K$ and, thus, $S_n$ converges to $\infty$ we can deduce from the tightness of $S_n^G$ given the data that $\phi_{n,\alpha}^G$ is consistent as well, compare to Theorem 7 of \cite{JanssenPauls2003}. 

Introduce $y_j=\kappa_j(1-G_j)(1-F_j)$ $(j=1,2)$ and $y=y_2+y_2$. A direct consequence of the extended Glivenko-Cantelli is $\sup\{|Y_{j}(t)/n-y_j(t)|:t\in[0,\infty)\}\to 0$ and, thus, $\sup\{|Y(t)/n-y(t)|:t\in[0,\infty)\}\to 0$, both in probability. The same argumentation leads to $\sup\{|N_{j}(t)/n-L_j(t)|:t\in[0,\infty)\}\to 0$ in probability, where  $L_j(t)=\kappa_j\P(X_{j,1}\leq t, \Delta_{j,1}=1)=\kappa_j\int_{[0,t]}(1-G_j)\,\mathrm{ d }F_j$ $(t\geq 0; \,j=1,2)$.  Moreover,
\begin{align*}
\widehat F(x) =1- \exp\Bigl( -\int_0^x \log(1-Y^{-1})^n \,\mathrm{ d }\frac{N}{n}\Bigr) \to 1- \exp\Bigl( -\int_0^x (1-y^{-1})\,\mathrm{ d }L\Bigr) \equiv F(x)\quad (x<\tau).
\end{align*}
By the monotonicity of the limit and standard subsequence arguments we can deduce $\sup\{|\widehat F(t)-F(t)|:t\in[0,x]\}\to 0$ in probability for every $x<\tau$. For every subsequence it is easy to construct a sub-subsequence that all uniform convergence results hold simultaneously with probability one. Since the final results do not depend on the pre-chosen subsequence we may operate in the following only along such sub-subsequences and on an appropriate event set with probability one. Hence, given the data we can treat the counting processes $N_j,N, Y_j,Y$ and the Kaplan--Meier estimator $\widehat F$ as (nonrandom) functions such that the previously mentioned uniform convergence statements hold. We do so in the following by always keeping the data fixed.

Let $X_{(1)}\leq \ldots \leq X_{(n)}$ be the order statistics of the pooled sample $(X_{j,i})_{j,i}$. Furthermore, let $G_{(i)}$ be the multiplier corresponding to $X_{(i)}$. Obviously, $G_{(1)},\ldots,G_{(n)}$ are still independent and identical distributed with the same distribution as $G_{1,1}$. Define
\begin{align*}
\xi_{n,i}^{(r)}= \Bigl(\frac{n}{n_1n_2}\Bigr)^{1/2} G_{(i)} w_r(\widehat F [X_{(i)}]) \frac{Y_1Y_2}{Y}(X_{(i)}) \Bigl[ \frac{\Delta N_1(X_{(i)})}{Y_1(X_{(i)})} - \frac{\Delta N_2(X_{(i)})}{Y_2(X_{(i)})}\Bigr]\quad (1\leq r \leq m;\,1\leq i \leq n).
\end{align*}
Thus we can represent the bootstrap version of ${\bf T}_n$ as ${\bf T}_n^G=\sum_{i=1}^n {\pmb\xi}_{n,i}$, where ${\pmb\xi}_{n,i}=(\xi_{n,i}^{(1)},\ldots,\xi_{n,n}^{(m)})^T$ $(1\le i \leq n)$. Using the multivariate Lindeberg-Feller theorem we can prove that it converges (given the data) in distribution to a multivariate normal distributed random vector ${\bf T}_0\sim N({\bf 0},{\bf \Sigma}_0)$ with covariance matrix ${\bf \Sigma}_0$ given by its entries
\begin{align*}
({\bf \Sigma}_0)_{r,s}= (\kappa_1\kappa_2)^{-1} \int (w_r\circ F) (w_s\circ F) \frac{y_1^2y_2^2}{y^2} (y_1^{-2}{\mathrm dL_1} + y_2^{-2}\mathrm dL_2) \quad(1\leq r,s\leq m).
\end{align*}
In particular, the covariance matrices ${\bf \Sigma}_0$ and ${\bf \Sigma}$ coincide under the null $\ho:A_1=A_2$ due to $\int_0^t y_j^{-1} \,\mathrm{ d }L_j=A_j(t)$. 

To accept the stated convergence first note that $\Delta N_1(X_{(i)})\Delta N_2(X_{(i)})$ equals $0$. With this in mind it is easy to check that
\begin{align*}
\sum_{i=1}^n \E(\xi_{n,i}^{(r)}\xi_{n,i}^{(s)}) = \frac{n}{n_1n_2} \int (w_r\circ \widehat F) (w_s\circ \widehat F) \frac{Y_1^2Y_2^2}{Y^2} \Bigl(\frac{\mathrm dN_1}{Y_1^2} + \frac{\mathrm dN_2}{Y_2^2}\Bigr)\to ({\bf \Sigma}_0)_{r,s}\quad(r,s=1,\ldots,m).
\end{align*}
It remains to verify the Lindeberg condition for each component. Since $|\xi_{n,i}^{(r)}|\leq \eta n^{-1/2}|G_i|$ for some $\eta>0$ independent of $i,r,n$ we obtain for every $\varepsilon>0$
\begin{align*}
\sum_{i=1}^n \E([\xi_{n,i}^{(r)}]^2\mathbf{1}\{|\xi_{n,i}^{(r)}|\geq\varepsilon\}) \leq \eta  \E(G_{(1)}^2 \mathbf{1}\{|G_{(1)}|\geq n^{1/2}\eta^{-1} \varepsilon\})\to 0
\end{align*}
by Lebesgue's dominated theorem and $\E(G_{(1)}^2)=1$. This proves that ${\bf T}_n^G$ converges in distribution to ${\bf T}_0\sim N({\bf 0},{\bf \Sigma}_0)$. Thus, it remains to analyze 
$\widehat{\bf \Sigma}^G_n$. To this end we introduce ${\bf \Sigma}_1\in\R^{m\times m}$ with entries
\begin{align*}
({\bf \Sigma}_1)_{r,s}&= (\kappa_1\kappa_2)^{-1} \int (w_r\circ F) (w_s\circ F) \frac{y_1y_2}{y} \frac{\mathrm d (L_1+L_2)}{y}\quad(r,s=1,\ldots,m).
\end{align*}
Under $\ho: F_1=F_2$ we have $\int y^{-1} \,\mathrm{ d }(L_1+L_2)= \int (1-F_1)^{-1} \,\mathrm{ d } F_1$ and, hence, ${\bf \Sigma}_1$ and ${\bf \Sigma}$ coincide under $\ho$. Moreover, it can easily checked that 
\begin{align*}
\E(({\bf\widehat\Sigma}^G_n)_{r,s}) = ({\bf\widehat\Sigma}_n)_{r,s} \to ({\bf \Sigma}_1)_{r,s}.
\end{align*}
Since $(a_1+\ldots+a_m)^2\leq m^2 (a_1^2+\ldots+a_m^2)$ for all $a_i\in\R$ and $w$ is bounded by $M$, say, we obtain
\begin{align*}
\E(({\bf\widehat\Sigma}^G_n)_{r,s}^2) \leq  m^2M^4\sum_{i=1}^n \E(G_{(i)}^4) \frac{n^2}{n_1^2n_2^2} = m^2M^4 \E(G_{(1)}^4) \frac{n^3}{n_1^2n_2^2}\to 0.
\end{align*}
Consequently, ${\bf\widehat\Sigma}^G_n\to {\bf \Sigma}_1$ in probability follows from Chebychev's inequality. 
Similar to the argumentation in Theorem~1 for ${\bf \Sigma}$, it follows that ${\bf \Sigma}_1$ is non-singular. Thus,  $\widehat{\bf \Sigma}_{n,J}^-$ converges to ${\bf \Sigma}_{1,J}^{-1}$ for every nonempty subset $J\subset\{1,\ldots,m\}$. Altogether we can conclude from the continuous mapping theorem that $S_n^G$ converges in distribution to the real-valued random variable 
\begin{align*}
S_0=\max\{0,{\bf T}^T_{0,J}{\bf \Sigma}_{1,J}^{-1}{\bf T}_{0,J}: \emptyset\neq J \subset \{1,\ldots,m\};\,{\bf \Sigma}_{1,J}^{-1}{\bf T}_{0,J}\geq 0\},
\end{align*}
where the distribution of $S_0$ equals the one of $S$ under $\ho$ Together with the arguments from the onset this completes the proof.
\qed

\section{Additional simulation results}
As stated in the main paper, we repeated the simulations under the null with $\widehat{\bf \Sigma}^G_n$ replaced by the empirical covariance matrix of ${\bf T}_n^G$ obtained in the $N_{boot}=1,000$ bootstrap Monte-Carlo runs (separately for each simulation run). The resulting empirical sizes are displayed in Table \ref{table:null_emp+var}. 
\begin{table}[ht]
	\caption{Empirical sizes in $\%$ (nominal level $5\%$) for three different wild bootstrap versions of our test using the empirical variance of ${\bf T}_n^G$ for studentization. The survival and censoring times were simulated by exponential distributions. }\label{table:null_emp+var}
	\centering
	\begin{tabular}{ccccc}
		\hline
		$(n_1,n_2)$ & censoring in \% & Normal & Poisson & Rademacher  \\ 
		\hline \hline
		& (10,30) & 8.7 & 6.6 & 8.1 \\ 
		(20,30)   & (15,15) & 10.5 & 8.3 & 9.5 \\ 
		& (30,30) & 10.7 & 7.8 & 9.9 \\
		\hline \hline
		& (10,30) & 6.8 & 4.9 & 5.8 \\ 
		(25,25)   & (15,15) & 8.7 & 6.3 & 7.3 \\  
		& (30,30) & 9.2 & 6.7 & 8.2 \\
		\hline\hline
		& (10,30) & 9.1 & 6.6 & 7.8 \\
		(30,70)   & (15,15) & 9.9 & 7.5 & 8.4 \\   
		& (30,30) & 10.8 & 8.3 & 9.4 \\ 
		\hline \hline
		& (10,30) & 7.2 & 4.9 & 5.8 \\ 
		(50,50)   & (15,15) & 7.9 & 5.7 & 6.4 \\ 
		& (30,30) & 7.8 & 5.4 & 6.2 \\ 
	\end{tabular}
\end{table}

Compared to the results from in Table \ref{table:null} in the main paper, the test based on Poisson multipliers now performs much better. It is still quite liberal but not as pronounced as with $\widehat{\bf \Sigma}^G_n$. Contrary, the tests corresponding to the Rademacher and normal multipliers are now more liberal.  Consequently, studentization by $\widehat{\bf\Sigma}_n^G$ should be preferred for them. Thus, leading to the interesting conjecture that the choice of studentization (empirical or based on squared multipliers) also depends on the choice of wild bootstrap multipliers. 
Anyhow, taking both results together, the Rademacher weights with studentization obtained via $\widehat{\bf\Sigma}_n^G$ showed the best results in our situation.

\section{GUI of \textit{mdir.onesided}}
\begin{figure}[h!] 
	\begin{center}	
		\includegraphics[width=0.3\textwidth]{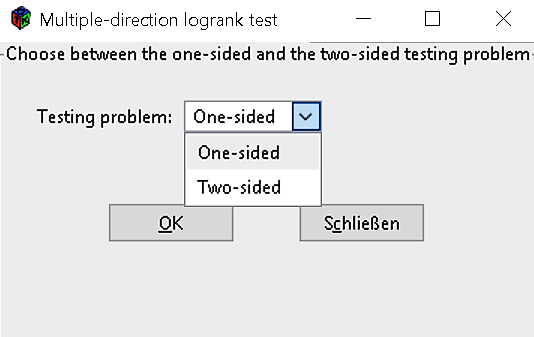}
	\end{center}
	\caption[h]{GUI: Decision between the one-sided and two-sided test.} \label{fig:GUI_one-twosided}
\end{figure}

\begin{figure}[h!] 
	\begin{center}	
		\includegraphics[width=0.5\textwidth]{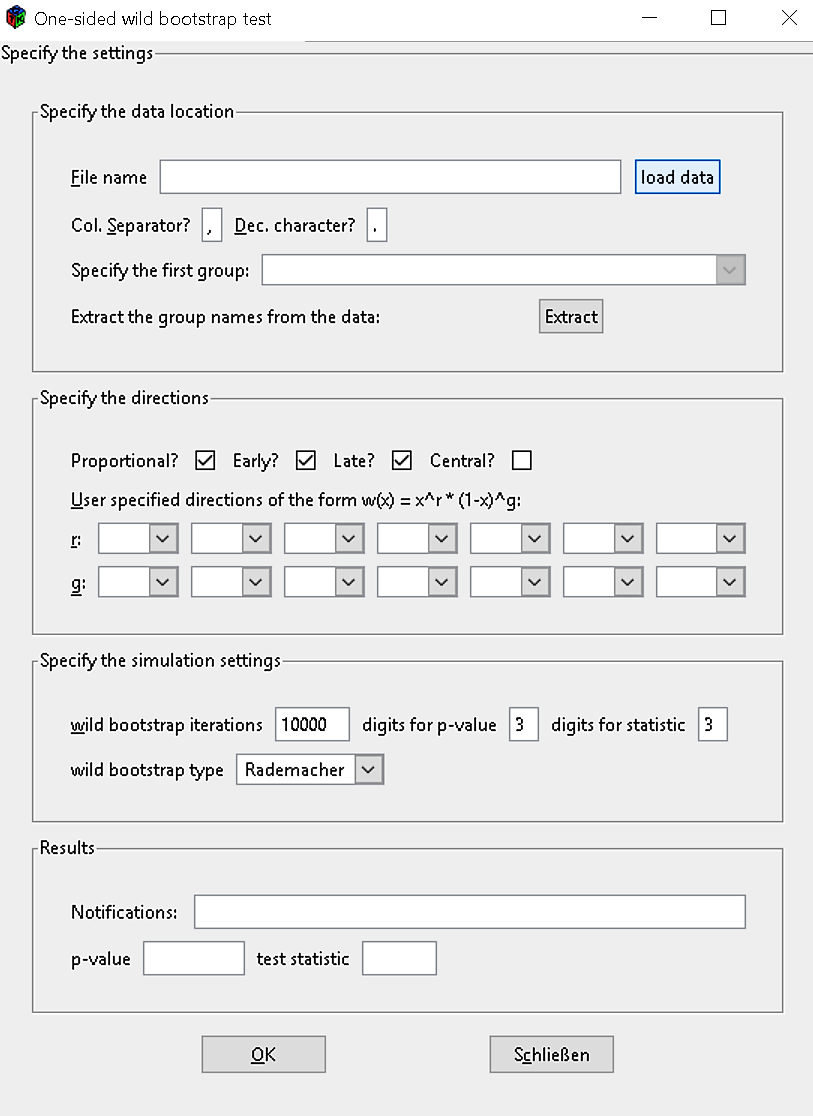}
	\end{center}
	\caption[h]{GUI: Empty window for the one-sided test. } \label{fig:GUI_empty}
\end{figure}

\begin{figure}[h!] 
	\begin{center}	
		\includegraphics[width=0.5\textwidth]{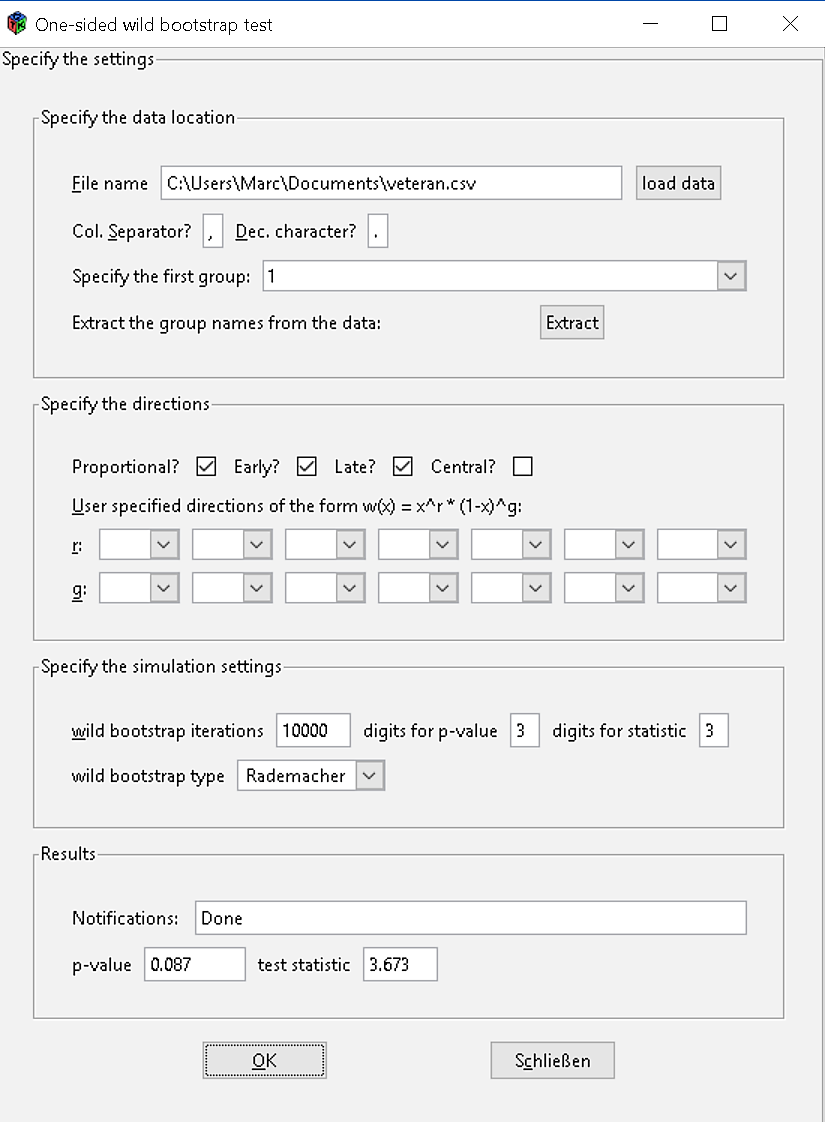} 	
	\end{center}
	\caption[h]{GUI: Application of our test with the default settings to the veteran data set. } \label{fig:GUI_default}
\end{figure}

\begin{figure}[h!] 
	\begin{center}	
		\includegraphics[width=0.5\textwidth]{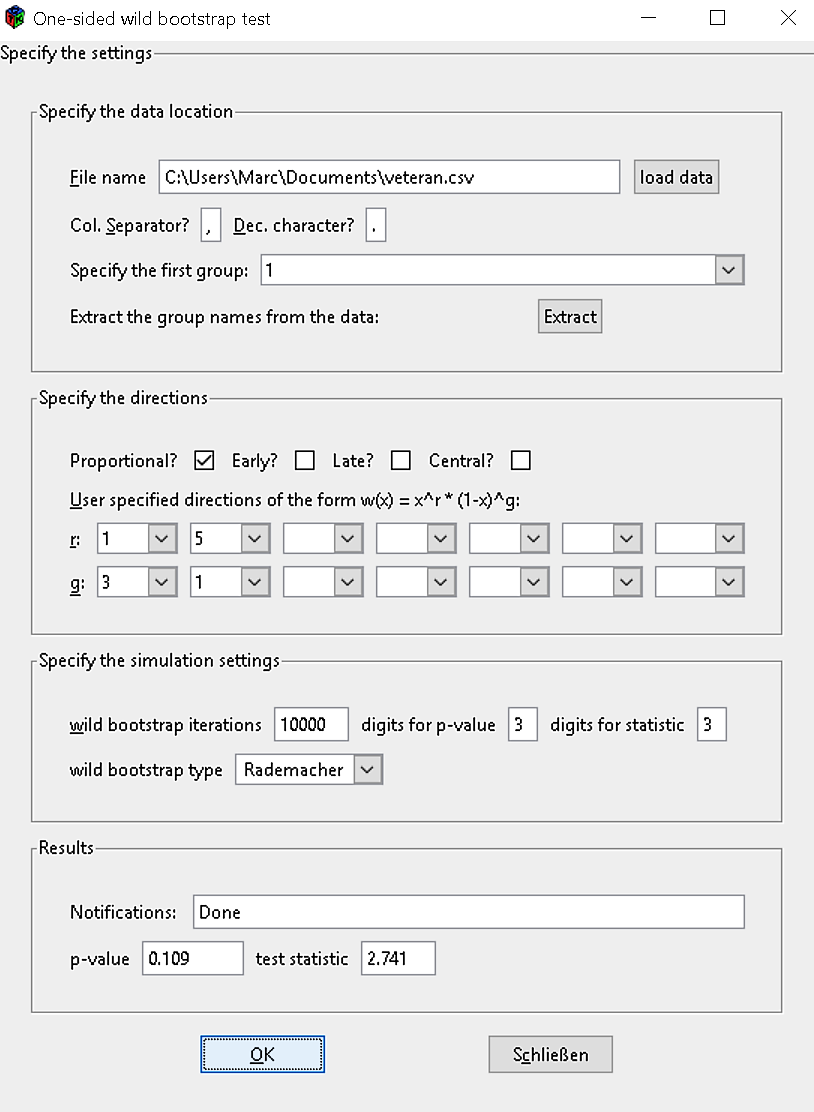}
	\end{center}
	\caption[h]{GUI: Application of our test with user specified directions to the veteran data set. } \label{fig:GUI_user-specified}
\end{figure}

For a simple use of our \textsc{R} function \textit{mdir.onesided} we implemented a graphical user interface (GUI). Note that our GUI only works
when the \textsc{R} package \textbf{RGtk2} is installed. After loading the \textsc{R} package \textbf{mdir.logrank} the command \textit{calculateGUI()} opens the first window (Figure \ref{fig:GUI_one-twosided}), where the user can choose between the one-sided and the two-sided test. 
After choosing \textit{One-sided} and clicking on the OK button the window for the one-sided test appears (Figure \ref{fig:GUI_empty}). 
The data set location can be specified directly or by the \textit{load data} button. To ensure that the data set is loaded correctly the user should check which symbols/characters are used for column separation and decimals, respectively, and indicate both in the two boxes intended for that. Keep in mind that rows consisting of the survival times, the group and censoring status need to be included in the chosen data set and need to be named by \textit{time}, \textit{group} and \textit{event}, respectively.  Via the \textit{Extract} button the groups' coding or the groups' names are extracted from the data set. The user then can choose the first group, i.e., the group corresponding to $F_1$ in our notation. Now, clicking the OK button starts the calculation for the default settings, i.e., the test based on the weights/directions $w_1,w_2,w_3$ which we used for our simulations in the main paper. For demonstration we load the veteran data set consisting of all tumor types. As already done in Section \ref{sec:realdata}, we test for the superiority of the experimental chemotherapy and, thus, choose the group coded by $1$ to be the first group. A screenshot for this set-up with the default settings is displayed in Figure \ref{fig:GUI_default}. Here, the slightly different $p$-value compared to Section~5 in the main paper ($0.087$ instead of $0.086$) is due to a different seed and corresponding Monte Carlo error. 
But the user can also change the wild bootstrap approach, the bootstrap runs and specify other weights of the form
\begin{align*}
w_{r,g}(x)=x^r(1-x)^g
\end{align*}
with exponents $r,g\in\{1,\ldots,20\}$. Note, that $1,000$ bootstrap runs with three chosen weights are calculated in a few seconds while $10,000$ runs need less than $1$ minute. Moreover, check boxes are provided for a quicker selection of the four common hazard directions:  proportional ($w_{0,0}$), early ($w_{0,4}$), late ($w_{4,0}$) and central ($w_{1,1}$). In Figure \ref{fig:GUI_user-specified} a further screenshot for the veteran data is given, where now the weights $w_{0,0}$, $w_{1,3}$ and $w_{5,1}$ are considered. 

\end{document}